\documentclass[a4paper]{amsart}
\usepackage{amssymb,amsrefs,mathrsfs,lineno,fancyvrb,tikz}
\usepackage{hyperref}

\newcommand\N{{\mathbb N}}
\newcommand\Z{{\mathbb Z}}
\newcommand\Q{{\mathbb Q}}
\newtheorem*{thm}{Theorem}

\DeclareMathOperator\Aut{Aut}
\DeclareMathOperator\Out{Out}
\DeclareMathOperator\sym{Sym}

\begin{document}
\title{The rational homology of the outer automorphism group of \boldmath $F_7$}
\author{Laurent Bartholdi}
\address{\'Ecole Normale Sup\'erieure, Paris\and Georg-August-Universit\"at zu G\"ottingen}
\email{laurent.bartholdi@gmail.com}
\date{18 January 2016}
\thanks{Partially supported by ANR grant ANR-14-ACHN-0018-01}
\begin{abstract}
  We compute the homology groups $H_*(\Out(F_7);\Q)$ of the outer
  automorphism group of the free group of rank $7$.

  We produce in this manner the first rational homology classes of
  $\Out(F_n)$ that are neither constant ($*=0$) nor Morita classes
  ($*=2n-4$).
\end{abstract}
\maketitle

\section{Introduction}
The homology groups $H_k(\Out(F_n);\Q)$ are intriguing objects. On the
one hand, they are known to ``stably vanish'', i.e. for all $n\in\N$
we have $H_k(\Out(F_n);\Q)=0$ as soon as $k$ is large
enough~\cite{galatius:stableautfn}. Hatcher and Vogtmann prove that
the natural maps $H_k\Out(F_n)\to H_k\Aut(F_n)$ and
$H_k\Aut(F_n)\to H_k\Aut(F_{n+1})$ are isomorphisms for $n\ge2k+2$
respectively $n\ge2k+4$,
see~\cites{hatcher-vogtmann:stability,hatcher-vogtmann-wahl:stability}. On
the other hand, $H_k(\Out(F_n);\Q)=0$ for $k>2n-3$, since $\Out(F_n)$
acts geometrically on a contractible space (the ``spine of outer
space'', see~\cite{culler-vogtmann:outerspace}) of dimension
$2n-3$. Combining these results, the only $k\ge1$ for which
$H_k(\Out(F_n);\Q)$ could possibly be non-zero are in the range
$\frac n2-2<k\le 2n-3$. Morita conjectures
in~\cite{morita:survey}*{page~390} that $H_{2n-3}(\Out(F_n);\Q)$
always vanishes; this would improve the upper bound to $k=2n-4$, and
$H_{2n-4}(\Out(F_n);\Q)$ is also conjectured to be non-trivial.

We shall see that the first conjecture does not hold.  Indeed, the
first few values of $H_k(\Out(F_n);\Q)$ may be computed by a
combination of human and computer work, and yield
\[\begin{array}{l|cccccccccccc}
n \backslash k & 0 & 1 & 2 & 3 & 4 & 5 & 6 & 7 & 8 & 9 & \hbox to 0pt{$10$}\phantom{1} & \hbox to 0pt{$11$}\phantom{1}\\ \hline
2 & 1 & 0\\
3 & 1 & 0 & 0 & 0\\
4 & 1 & 0 & 0 & 0 & 1 & 0\\
5 & 1 & 0 & 0 & 0 & 0 & 0 & 0 & 0\\
6 & 1 & 0 & 0 & 0 & 0 & 0 & 0 & 0 & 1 & 0\\
7 & 1 & 0 & 0 & 0 & 0 & 0 & 0 & 0 & 1 & 0 & 0 & 1\\
\end{array}\]

The values for $n\le6$ were computed by Ohashi
in~\cite{ohashi:OutF6}. They reveal that, for $n\le6$, only the
constant class ($k=0$) and the Morita classes $k=2n-4$ yield
non-trivial homology. The values for $n=7$ are the object of this
Note, and reveal that the picture changes radically:
\begin{thm}
  The non-trivial homology groups $H_k(\Out(F_7);\Q)$ occur for
  $k\in\{0,8,11\}$ and are all $1$-dimensional.
\end{thm}

Previously, only the rational Euler characteristic
$\chi_\Q(\Out(F_7))=\sum(-1)^k\dim H_k(\Out(F_7);\Q)$ was
known~\cite{morita-sakasai-suzuki:computations}, and shown to be
$1$. These authors computed in fact the rational Euler characteristics
up to $n=11$ in that paper and the
sequel~\cite{morita-sakasai-suzuki:OutF11}.

\section{Methods}
We make fundamental use of a construction of
Kontsevich~\cite{kontsevich:formalsymplectic}, explained
in~\cite{conant-vogtmann:kontsevich}. We follow the simplified description
from~\cite{ohashi:OutF6}.

Let $F_n$ denote the free group of rank $n$. This parameter $n$ is fixed
once and for all, and will in fact be omitted from the notation as often as
possible. An \emph{admissible graph of rank $n$} is a graph $G$ that is
$2$-connected ($G$ remains connected even after an arbitrary edge is
removed), without loops, with fundamental group isomorphic to $F_n$, and
without vertex of valency $\le2$. Its \emph{degree} is $\deg(G):=\sum_{v\in
  V(G)}(\deg(v)-3)$. In particular, $G$ has $2n-2-\deg(G)$ vertices and
$3n-3-\deg(G)$ edges, and is trivalent if and only if $\deg(G)=0$. If
$\Phi$ is a collection of edges in a graph $G$, we denote by $G/\Phi$ the
graph quotient, obtained by contracting all edges in $\Phi$ to points.

A \emph{forested graph} is a pair $(G,\Phi)$ with $\Phi$ an oriented forest
in $G$, namely an ordered collection of edges that do not form any
cycle. We note that the symmetric group $\sym(k)$ acts on the set of
forested graphs whose forest contains $k$ edges, by permuting the forest's
edges.

For $k\in\N$, let $C_k$ denote the $\Q$-vector space spanned by
isomorphism classes of forested graphs of rank $n$ with a forest of
size $k$, subject to the relation
\[(G,\pi\Phi)=(-1)^\pi(G,\Phi)\text{ for all }\pi\in\sym(k).
\]
Note, in particular, that if $(G,\Phi)\sim(G,\pi\Phi)$ for an odd
permutation $\pi$ then $(G,\Phi)=0$ in $C_k$. These spaces $(C_*)$
form a chain complex for the differential
$\partial=\partial_C-\partial_R$, defined respectively on
$(G,\Phi)=(G,\{e_1,\dots,e_p\})$ by
\begin{align*}
  \partial_C(G,\Phi)&=\sum_{i=1}^p(-1)^i(G/e_i,\Phi\setminus\{e_i\}),\\
  \partial_R(G,\Phi)&=\sum_{i=1}^p(-1)^i(G,\Phi\setminus\{e_i\}),
\end{align*}
and the homology of $(C_*,\partial)$ is $H_*(\Out(F_n);\Q)$.

The spaces $C_k$ may be filtered by degree: let $F_p C_k$ denote the
subspace spanned by forested graphs $(G,\Phi)$ with $\deg(G/\Phi)\le p$. The
differentials satisfy respectively
\[\partial_C(F_p C_k)\subseteq F_p C_{k-1},\qquad\partial_R(F_p C_k)\subseteq F_{p-1}C_{k-1}.\]
A spectral sequence argument gives
\begin{equation}\label{eq:ss}
  H_p(\Out(F_n);\Q)=E^2_{p,0}=\frac{\ker(\partial_C|F_p C_p)\cap
    \ker(\partial_R|F_p C_p)}{\partial_R(\ker(\partial_C|F_{p+1}C_{p+1}))}.
\end{equation}
Note that if $(G,\Phi)\in F_p C_p$ then $G$ is trivalent.  We compute
explicitly bases for the vector spaces $F_p C_p$, and matrices for the
differentials $\partial_C,\partial_R$, to prove the theorem.

\section{Implementation}
We follow for $n=7$ the procedure sketched
in~\cite{ohashi:OutF6}. Using the software program
\texttt{nauty}~\cite{mckay-piperno:nauty}, we enumerate all trivalent
graphs of rank $n$ and vertex valencies $\ge3$. The libraries in
\texttt{nauty} produce a canonical ordering of a graph, and compute
generators for its automorphism group. We then weed out the
non-$2$-connected ones.

For given $p\in\N$, we then enumerate all $p$-element oriented forests in
these graphs, and weed out those that admit an odd symmetry. These are
stored as a basis for $F_p C_p$. Let $a_p$ denote the dimension of $F_p
C_p$.

For $(G,\Phi)$ a basis vector in $F_p C_p$, the forested graphs that
appear as summands in $\partial_C(G,\Phi)$ and $\partial_R(G,\Phi)$
are numbered and stored in a hash table as they occur, and the
matrices $\partial_C$ and $\partial_R$ are computed as sparse matrices
with $a_p$ columns.

The nullspace $\ker(\partial_C|F_p C_p)$ is then computed: let $b_p$
denote its dimension; then the nullspace is stored as a sparse
$(a_p\times b_p)$-matrix $N_p$. The computation is greatly aided by
the fact that $\partial_C$ is a block matrix, whose row and column
blocks are spanned by $\{(G,\Phi):G/\Phi=G_0\}$ for all choices of the
fully contracted graph $G_0$. The matrices $N_p$ are computed using
the linear algebra library \texttt{linbox}~\cite{linbox}, which
provides exact linear algebra over $\Q$ and finite fields.

Finally, the rank $c_p$ of $\partial_R\circ N_p$ is computed, again using
\texttt{linbox}.  By~\eqref{eq:ss}, we have
\[\dim H_p(\Out(F_n);\Q)= b_p - c_p - c_{p+1}.\]

For memory reasons (the computational requirements reached 200GB of
RAM at its peak), some of these ranks were computed modulo a large
prime ($65521$ and $65519$ were used in two independent runs).

Computing modulo a prime can only reduce the rank; so that the values $c_p$
we obtained are underestimates of the actual ranks of $\partial_R\circ
N_p$. However, we also know \emph{a priori} that $b_p - c_p - c_{p+1}\ge0$
since it is the dimension of a vector space; and none of the $c_p$ we
computed can be increased without at the same time causing a homology
dimension to become negative, so our reduction modulo a prime is legal.

For information, the parameters $a_p,b_p,c_p$ for $n=7$ are as follows:
\[\begin{array}{r|llllllllllll}
  p & 0 & 1 & 2 & 3 & 4 & 5 & 6 & 7 & 8 & 9 & 10 & 11\\\hline
  a_p & 365 & 3712 & 23227 & {\approx}105k & {\approx}348k & {\approx}854k & {\approx}1.6m & {\approx}2.3m & {\approx}2.6m & {\approx}2.1m & {\approx}1.2m & {\approx}376k\\
  b_p & 365 & 1784 & 5642 & 14766 & 28739 & 39033 & 38113 & 28588 & 16741 &
  6931& 1682 & 179\\
  c_p & 0 & 364 & 1420 & 4222 & 10544 & 18195 & 20838 & 17275 & 11313 &
  5427 & 1504 & 178
\end{array}\]

The largest single matrix operations that had to be performed were
computing the nullspace of a $2038511\times536647$ matrix (16 CPU hours)
and the rank modulo $65519$ of a (less sparse) $1355531\times16741$ matrix
(10 CPU hours).

The source files used for the computations are available as
supplemental material. Compilation requires \texttt{g++} version 4.7
or later, a functional \texttt{linbox} library, available from the
site \texttt{http://www.linalg.org}, as well as the \texttt{nauty}
program suite, available from the site
\texttt{http://pallini.di.uniroma1.it}. It may also be directly
downloaded and installed by typing `\texttt{make nauty25r9}' in the
directory in which the sources were downloaded. Beware that the
calculations required for $n=7$ are prohibitive for most desktop
computers.

\section*{Conclusion}
Computing the dimensions of the homology groups is only the first step
in understanding them; much more interesting would be to know
visually, or graph-theoretically, where these non-trivial classes come
from.

It seems almost hopeless to describe, via computer experiments, the
non-trivial class in degree $8$. It may be possible, however, to
arrive at a reasonable understanding of the non-trivial class in
degree $11$.

This class may be interpreted as a linear combination $w$ of trivalent
graphs on $12$ vertices, each marked with an oriented spanning
forest. There are $376365$ such forested graphs that do not admit an
odd symmetry. The class $w\in\Q^{376365}$ is an $\Z$-linear
combination of $70398$ different forested graphs, with coefficients in
$\{\pm1,\dots,\pm16\}$. For example, eleven graphs occur with
coefficient $\pm13$; four of them have indices
$25273, 53069, 53239, 53610$ respectively, and are, with the spanning
tree in bold,
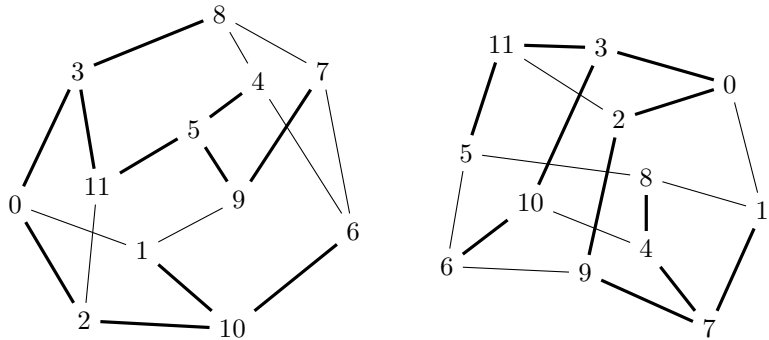
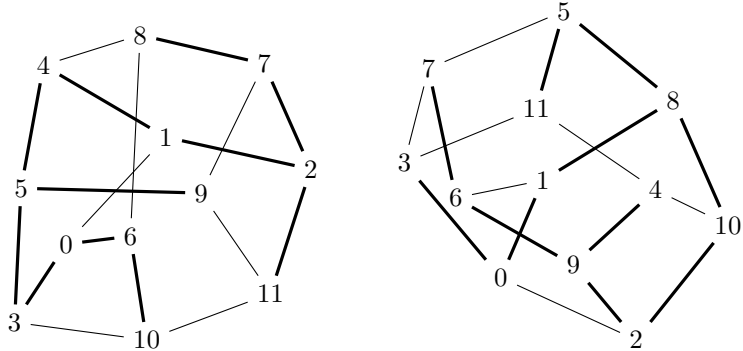
\begin{center}
\begin{tikzpicture}[scale=0.6]
  \node (11) at (-60bp,-6bp) {11};
  \node (10) at (24bp,-95bp) {10};
  \node (1) at (-32bp,-46bp) {1};
  \node (0) at (-111bp,-17bp) {0};
  \node (3) at (-72bp,66bp) {3};
  \node (2) at (-68bp,-90bp) {2};
  \node (5) at (0bp,30bp) {5};
  \node (4) at (40bp,60bp) {4};
  \node (7) at (80bp,66bp) {7};
  \node (6) at (99bp,-34bp) {6};
  \node (9) at (28bp,-14bp) {9};
  \node (8) at (16bp,102bp) {8};

  \draw[very thick] (0) -- (3);
  \draw[very thick] (3) -- (11);
  \draw[very thick] (4) -- (5);
  \draw[very thick] (0) -- (2);
  \draw[very thick] (7) -- (9);
  \draw[very thick] (3) -- (8);
  \draw[thin] (4) -- (8);
  \draw[thin] (6) -- (7);
  \draw[thin] (4) -- (6);
  \draw[thin] (7) -- (8);
  \draw[very thick] (5) -- (11);
  \draw[very thick] (1) -- (10);
  \draw[very thick] (6) -- (10);
  \draw[thin] (2) -- (11);
  \draw[thin] (0) -- (1);
  \draw[thin] (1) -- (9);
  \draw[very thick] (2) -- (10);
  \draw[very thick] (5) -- (9);
\end{tikzpicture}
\qquad\begin{tikzpicture}[scale=0.6]
  \node (11) at (-60bp,93bp) {11};
  \node (10) at (-42bp,-6bp) {10};
  \node (1) at (102bp,-11bp) {1};
  \node (0) at (82bp,68bp) {0};
  \node (3) at (2bp,91bp) {3};
  \node (2) at (13bp,46bp) {2};
  \node (5) at (-82bp,25bp) {5};
  \node (4) at (30bp,-35bp) {4};
  \node (7) at (69bp,-85bp) {7};
  \node (6) at (-94bp,-46bp) {6};
  \node (9) at (-8bp,-50bp) {9};
  \node (8) at (30bp,10bp) {8};
  \draw[very thick] (0) -- (3);
  \draw[very thick] (3) -- (11);
  \draw[thin] (5) -- (6);
  \draw[very thick] (0) -- (2);
  \draw[thin] (5) -- (8);
  \draw[very thick] (7) -- (9);
  \draw[very thick] (4) -- (8);
  \draw[thin] (2) -- (11);
  \draw[thin] (4) -- (10);
  \draw[very thick] (5) -- (11);
  \draw[very thick] (2) -- (9);
  \draw[very thick] (1) -- (7);
  \draw[very thick] (6) -- (10);
  \draw[thin] (1) -- (8);
  \draw[thin] (0) -- (1);
  \draw[thin] (6) -- (9);
  \draw[very thick] (4) -- (7);
  \draw[very thick] (3) -- (10);
\end{tikzpicture}\\
\begin{tikzpicture}[scale=0.6]
  \node (11) at (77bp,-66bp) {11};
  \node (10) at (0bp,-95bp) {10};
  \node (1) at (12bp,32bp) {1};
  \node (0) at (-50bp,-35bp) {0};
  \node (3) at (-82bp,-84bp) {3};
  \node (2) at (102bp,12bp) {2};
  \node (5) at (-78bp,0bp) {5};
  \node (4) at (-64bp,77bp) {4};
  \node (7) at (73bp,78bp) {7};
  \node (6) at (-10bp,-30bp) {6};
  \node (9) at (34bp,-3bp) {9};
  \node (8) at (-4bp,96bp) {8};
  \draw[very thick] (0) -- (3);
  \draw[thin] (4) -- (8);
  \draw[very thick] (4) -- (5);
  \draw[very thick] (2) -- (7);
  \draw[thin] (7) -- (9);
  \draw[very thick] (5) -- (9);
  \draw[very thick] (6) -- (10);
  \draw[very thick] (7) -- (8);
  \draw[thin] (10) -- (11);
  \draw[very thick] (1) -- (4);
  \draw[very thick] (1) -- (2);
  \draw[very thick] (0) -- (6);
  \draw[very thick] (2) -- (11);
  \draw[thin] (0) -- (1);
  \draw[thin] (9) -- (11);
  \draw[very thick] (3) -- (5);
  \draw[thin] (3) -- (10);
  \draw[thin] (6) -- (8);
\end{tikzpicture}
\qquad\begin{tikzpicture}[scale=0.6]
  \node (11) at (-10bp,51bp) {11};
  \node (10) at (109bp,-22bp) {10};
  \node (1) at (-6bp,6bp) {1};
  \node (0) at (-32bp,-55bp) {0};
  \node (3) at (-92bp,17bp) {3};
  \node (2) at (52bp,-94bp) {2};
  \node (5) at (7bp,112bp) {5};
  \node (4) at (64bp,0bp) {4};
  \node (7) at (-77bp,76bp) {7};
  \node (6) at (-60bp,-5bp) {6};
  \node (9) at (13bp,-47bp) {9};
  \node (8) at (75bp,56bp) {8};
  \draw[very thick] (0) -- (3);
  \draw[thin] (3) -- (11);
  \draw[thin] (0) -- (2);
  \draw[very thick] (5) -- (8);
  \draw[thin] (5) -- (7);
  \draw[thin] (4) -- (11);
  \draw[thin] (4) -- (10);
  \draw[very thick] (6) -- (7);
  \draw[thin] (1) -- (6);
  \draw[very thick] (5) -- (11);
  \draw[very thick] (2) -- (9);
  \draw[thin] (3) -- (7);
  \draw[very thick] (8) -- (10);
  \draw[very thick] (1) -- (8);
  \draw[very thick] (0) -- (1);
  \draw[very thick] (2) -- (10);
  \draw[very thick] (6) -- (9);
  \draw[very thick] (4) -- (9);
\end{tikzpicture}
\end{center}

The coefficients of $w$, and corresponding graphs, are distributed as
ancillary material in the file \verb+w_cycle+, in format
`\verb+coefficient [edge1 edge2 ...]+', where each edge is
`\verb|x-y|' or `\verb|x+y|' to indicate whether the edge is absent or
present in the forest. Edges always satisfy $\verb+x+\le\verb+y+$, and
the forest is oriented so that its edges are lexicographically
ordered. Edges are numbered from $0$ while graphs are numbered from
$1$. There are no multiple edges.

\section*{Acknowledgments}
I am grateful to Alexander Berglund and Nathalie Wahl for having
organized a wonderful and stimulating workshop on automorphisms of
free groups in Copenhagen in October 2015, when this work began; to
Masaaki Suzuki, Andy Putman and Karen Vogtmann for very helpful
conversations that took place during this workshop; and to Jim Conant
for having checked the cycle $w$ (after finding a mistake in its
original signs) with an independent program.


\end{document}